\documentclass[11pt]{article}
\usepackage{amsfonts}
\usepackage{graphics}
\usepackage{indentfirst}
\usepackage{color}
\usepackage{cite}
\usepackage{latexsym}
\usepackage[paper=a4paper, left=2.1cm, right=2.1cm, top=2.2cm, bottom=1.5cm, headheight=5.5pt, footskip=0.8cm, footnotesep=0.8cm, centering, includefoot]{geometry}
\usepackage{amsmath}
\allowdisplaybreaks
\usepackage{amssymb}
\usepackage[dvips]{epsfig}
\usepackage{amscd}

\newtheorem{theorem}{Theorem}[section]
\newtheorem{remark}{Remark}[section]
\newtheorem{definition}{Definition}[section]
\newtheorem{lemma}{Lemma}[section]

\DeclareMathOperator{\divv}{div}

\makeatletter
\@addtoreset{equation}{section}
\makeatother
\makeatletter
\@addtoreset{equation}{section}
\makeatother

\title{{Singularity formation to the two-dimensional full compressible Navier-Stokes equations with zero heat conduction in a bounded domain}
\thanks{Supported by Chongqing
Research Program of Basic Research and Frontier Technology (No. cstc2018jcyjAX0049),  the Postdoctoral Science Foundation of Chongqing (No. xm2017015), and China Postdoctoral Science Foundation
(Nos. 2018T110936, 2017M610579).}
}

\author{Xin Zhong\thanks{School of Mathematics and Statistics, Southwest University, Chongqing 400715,
People's Republic of China ({\tt xzhong1014@amss.ac.cn}).
}
}
\date{ }

\begin{document}
\maketitle

\begin{abstract}
We consider the singularity formation of strong solutions to the two-dimensional full compressible Navier-Stokes equations with zero heat conduction in a bounded domain. It is shown that for the initial density allowing vacuum, the strong solution exists globally if the density and the pressure are bounded from above.
Critical Sobolev inequalities of logarithmic type play a crucial role in the proof.
\end{abstract}

Keywords: full compressible Navier-Stokes equations; zero heat conduction; blow-up criterion.

Math Subject Classification: 35Q30; 35B65

\section{Introduction}
Let $\Omega\subset\mathbb{R}^n\ (n=2,3)$ be a domain, the motion of a viscous, compressible, and heat conducting Navier-Stokes flow in $\Omega$ can be described by the full compressible Navier-Stokes equations
\begin{align}\label{1.1}
\begin{cases}
\rho_{t}+\divv(\rho\mathbf{u})=0,\\
(\rho\mathbf{u})_{t}+\divv(\rho\mathbf{u}\otimes\mathbf{u})
-\mu\Delta\mathbf{u}
-(\lambda+\mu)\nabla\divv\mathbf{u}+\nabla P=
\mathbf{0},\\
c_{\nu}[(\rho\theta)_{t}+\divv(\rho\mathbf{u}\theta)]
+P\divv\mathbf{u}-\kappa\Delta\theta
=2\mu|\mathfrak{D}(\mathbf{u})|^2+\lambda(\divv\mathbf{u})^2.
\end{cases}
\end{align}
Here, $t\geq0$ is the time, $x\in\Omega$ is the spatial coordinate, and the unknown functions $\rho, \mathbf{u}, P=R\rho\theta\ (R>0), \theta$ are the fluid density, velocity, pressure, and the absolute temperature respectively; $\mathfrak{D}(\mathbf{u})$ denotes the deformation tensor given by
\begin{equation*}
\mathfrak{D}(\mathbf{u})=\frac{1}{2}(\nabla\mathbf{u}+(\nabla\mathbf{u})^{tr}).
\end{equation*}
The constant viscosity coefficients $\mu$ and $\lambda$ satisfy the physical restrictions
\begin{equation}\label{1.2}
\mu>0,\ 2\mu+n\lambda\geq0.
\end{equation}
Positive constants $c_\nu$ and $\kappa$ are respectively the heat capacity, the ratio of the heat conductivity coefficient over the heat capacity.

The study on the theory of well-posedness of solutions to the Cauchy problem and the initial boundary value problem (IBVP) for the compressible Navier-Stokes equations has grown enormously in recent years,
refer to \cite{CK2006,H1995,H1997,HLX2012,HL2013,L1998,HL2018,
LX2013,LL2014,YZ2017,WZ2017,MN1980,MN1983} and references therein.
In particular, non-vacuum small perturbations of a uniform non-vacuum
constant state have been shown existing globally in time and remain smooth in any space dimensions \cite{MN1980,MN1983}, while for general data which may contain vacuum states, only weak solutions are shown to exist for the compressible Navier-Stokes system in multi-dimension with special equation of state as in \cite{FNP2001,L1998}, yet the uniqueness
and regularity of these weak solutions remain unknown.
Despite the surprising results on global well-posedness of the strong (or classical) solution to the multi-dimensional compressible Navier-Stokes system for initial data with small total energy but possible large oscillations and containing vacuum states \cite{HLX2012,HL2018,LX2013,YZ2017,WZ2017},
it is an outstanding challenging open problem to investigate the global well-posedness for general large strong solutions with vacuum.


Therefore, it is important to study the mechanism of blow-up and structure of possible singularities of strong (or classical) solutions to the compressible Navier-Stokes equations. The pioneering work can be traced
to Serrin's criterion \cite{S1962} on the Leray-Hopf weak solutions to the three-dimensional incompressible Navier-Stokes equations, which can be stated
that if a weak solution $\mathbf{u}$ satisfies
\begin{equation}\label{1.3}
\mathbf{u}\in L^s(0,T;L^r),\ \text{for}\ \frac2s+\frac3r=1,\ 3<r\leq\infty,
\end{equation}
then it is regular. For more information on the blow-up criteria of compressible Navier-Stokes equations, we refer to \cite{WZ2013,HLX20112,SWZ20111} and references therein.

Recently, there are several results on the blow-up criteria of strong (or classical) solutions to the full compressible Navier-Stokes equations \eqref{1.1}. Precisely, let $0<T^*<+\infty$ be the maximum time of existence of strong solutions.
For the 3D full compressible Navier-Stokes equations,
under the assumption
\begin{equation*}
7\mu>\lambda,
\end{equation*}
Fan-Jiang-Ou \cite{FJO2010} showed that
\begin{equation}\label{1.4}
\lim_{T\rightarrow T^{*}}\left(\|\nabla\mathbf{u}\|_{L^1(0,T;L^\infty)}+\|\theta\|_{L^\infty(0,T;L^\infty)}\right)
=\infty.
\end{equation}
It was improved later in \cite{WZ2013} to
\begin{equation}\label{1.04}
\lim_{T\rightarrow T^{*}}\left(\|\rho\|_{L^\infty(0,T;L^\infty)}+\|\theta\|_{L^\infty(0,T;L^\infty)}\right)
=\infty
\end{equation}
provided that
\begin{equation*}
3\mu>\lambda.
\end{equation*}
Under just the physical condition
\begin{equation*}
\mu>0,\ 2\mu+3\lambda\geq0,
\end{equation*}
Huang-Li-Wang \cite{HLW2013} established the following Serrin type criterion
\begin{equation}\label{1.5}
\lim_{T\rightarrow T^{*}}\left(\|\divv\mathbf{u}\|_{L^1(0,T;L^\infty)}
+\|\mathbf{u}\|_{L^s(0,T;L^r)}\right)=\infty, \ \text{for}\ \ \frac2s+\frac3r=1,\ 3<r\leq\infty.
\end{equation}
Later on, Huang-Li \cite{HL2013}  extended \eqref{1.5} with a refiner form
\begin{equation}\label{1.6}
\lim_{T\rightarrow T^{*}}\left(\|\rho\|_{L^\infty(0,T;L^\infty)}+\|\mathbf{u}\|_{L^s(0,T;L^r)}\right)
=\infty,\ \text{for}\ \frac2s+\frac3r\leq1,\ 3<r\leq\infty.
\end{equation}
On the other hand, for the IBVP of 2D full Navier-Stokes equations, Wang \cite{W2014} showed the formation of singularity must be caused by losing the bound of $\divv\mathbf{u}$. More precisely, she obtained
\begin{equation}\label{1.7}
\lim_{T\rightarrow T^{*}}\|\divv\mathbf{u}\|_{L^1(0,T;L^\infty)}=\infty.
\end{equation}

It should be noted that all the results mentioned above on the blow-up criteria of strong (or classical) solutions of full compressible Navier-Stokes equations are for $\kappa>0$. Recently, for the Cauchy problem and the IBVP of 3D full compressible Navier-Stokes equations \eqref{1.1} with $\kappa=0$,
under the assumption
\begin{equation}\label{1.07}
\mu>4\lambda,
\end{equation}
Huang-Xin \cite{HX2016} showed that
\begin{equation}\label{1.06}
\lim_{T\rightarrow T^{*}}\left(\|\rho\|_{L^\infty(0,T;L^\infty)}
+\|\theta\|_{L^\infty(0,T;L^\infty)}\right)
=\infty,
\end{equation}
while for the Cauchy problem of 2D case, Zhong \cite{Z2017} established the following criterion
\begin{align}\label{1.05}
\lim_{T\rightarrow T^{*}}\left(\|\rho\|_{L^{\infty}(0,T;L^\infty)}
+\|P\|_{L^{\infty}(0,T;L^\infty)}\right)=\infty.
\end{align}
It is worth mentioning that in a well-known paper \cite{X1998}, Xin considered full compressible Navier-Stokes equations with $\kappa=0$ in multidimensional space, starting with a compactly supported initial density. He first proved that if the support of the density grows sublinearly in time and if the entropy is bounded from below then the solution cannot exist for all time. One key ingredient in the proof is a differential inequality on some integral functional (see \cite[Proposition 2.1]{X1998} for details). As an application, any smooth solution to the full compressible Navier-Stokes equations for polytropic fluids in the absence of heat conduction will blow up in finite time if the initial density is compactly supported.
Recently, based on the key observation that if initially a positive mass
is surrounded by a bounded vacuum region, then the time evolution remains uniformly bounded for all time, Xin-Yan \cite{XY2013} improved the blow-up results in \cite{X1998} by removing the assumptions that the initial density has compact support and the smooth solution has finite energy, but the initial data only has an isolated mass group.
Thus it seems very difficult to study globally smooth solutions of full compressible Navier-Stokes equations without heat conductivity in
multi-dimension. These motivate us to study a blow-up criterion for the IBVP of 2D full Navier-Stokes equations with zero heat conduction.
In fact, this is the main aim of this paper.


Let $\Omega\subset\mathbb{R}^2$ be a bounded smooth domain, when $\kappa=0$, and without loss of generality, take $c_\nu=R=1$, then the system \eqref{1.1} can be written as
\begin{align}\label{1.10}
\begin{cases}
\rho_{t}+\divv(\rho\mathbf{u})=0,\\
(\rho\mathbf{u})_{t}+\divv(\rho\mathbf{u}\otimes\mathbf{u})
-\mu\Delta\mathbf{u}
-(\lambda+\mu)\nabla\divv\mathbf{u}+\nabla P=
\mathbf{0},\\
P_{t}+\divv(P\mathbf{u})
+P\divv\mathbf{u}
=2\mu|\mathfrak{D}(\mathbf{u})|^2+\lambda(\divv\mathbf{u})^2,
\end{cases}
\end{align}
and the constant viscosity coefficients $\mu$ and $\lambda$ satisfy
\begin{equation}\label{1.02}
\mu>0,\ \mu+\lambda\geq0.
\end{equation}
The present paper aims at giving a blow-up criterion of strong solutions to the system \eqref{1.10} with the initial condition
\begin{equation}\label{1.11}
(\rho,\rho\mathbf{u},P)(x,0)=(\rho_0,\rho_0\mathbf{u}_0,P_0)(x),\ \ x\in\Omega,
\end{equation}
and the boundary condition
\begin{equation}\label{1.12}
\mathbf{u}=\mathbf{0},\  \ x\in\partial\Omega.
\end{equation}

Before stating our main result, we first explain the notations and conventions used throughout this paper. Set
\begin{equation*}
\int \cdot dx\triangleq\int_{\Omega}\cdot dx.
\end{equation*}
For $1\leq p\leq\infty$ and integer $k\geq0$, the standard Sobolev spaces are denoted by:
\begin{equation*}
L^p=L^p(\Omega),\ W^{k,p}=W^{k,p}(\Omega), \ H^{k}=H^{k,2}(\Omega),\ H_0^1=\{u\in H^1:u|_{\partial\Omega}=0\}.
\end{equation*}
Now we define precisely what we mean by strong solutions to the problem \eqref{1.10}--\eqref{1.12}.
\begin{definition}[Strong solutions]\label{def1}
$(\rho,\mathbf{u},P)$ is called a strong solution to \eqref{1.10}--\eqref{1.12} in $\Omega\times(0,T)$, if for some $q>2$,
\begin{align*}
\begin{cases}
\rho\geq0,\ \rho\in C([0,T]; W^{1,q}),\ \rho_t\in C([0,T]; L^{q}), \\
\nabla\mathbf{u}\in L^{\infty}(0,T; H^1)\cap L^2(0,T;W^{1,q}),\\
\sqrt{\rho}\mathbf{u}, \ \sqrt{\rho}\dot{\mathbf{u}}\in L^{\infty}(0,T; L^2), \\
P\geq0,\ P\in C([0,T];W^{1,q}),\ P_t\in C([0,T];L^{q}),
\end{cases}
\end{align*}
and $(\rho,\mathbf{u},P)$ satisfies both \eqref{1.10} almost everywhere in $\Omega\times(0,T)$ and \eqref{1.11} almost everywhere in $\Omega$.
\end{definition}

Our main result reads as follows:
\begin{theorem}\label{thm1.1}
Assume that the initial data $(\rho_0\geq0, \mathbf{u}_0,P_0\geq0)$ satisfies for any given $q>2$,
\begin{equation}\label{2.02}
\rho_{0}\in W^{1,q},\ \mathbf{u}_{0}\in H_0^1\cap H^2,\ \ P_0\in W^{1,q},
\end{equation}
and the compatibility condition
\begin{equation}\label{A2}
-\mu\Delta\mathbf{u}_0-(\lambda+\mu)\nabla\divv\mathbf{u}_0+\nabla P_0=\sqrt{\rho_0}\mathbf{g}
\end{equation}
for some $\mathbf{g}\in L^2(\Omega)$.
Let $(\rho,\mathbf{u},P)$ be a strong solution to the problem \eqref{1.10}--\eqref{1.12}. If $T^{*}<\infty$ is the maximal time of existence for that solution,
then we have
\begin{align}\label{B}
\lim_{T\rightarrow T^{*}}\left(\|\rho\|_{L^{\infty}(0,T;L^\infty)}
+\|P\|_{L^{\infty}(0,T;L^\infty)}\right)=\infty.
\end{align}
\end{theorem}

Several remarks are in order.
\begin{remark}\label{re1.1}
The local existence of a strong solution with initial data as in Theorem \ref{thm1.1} can be established by similar arguments as in \cite{CK2006}. Thus, the maximal time $T^{*}$ is well-defined.
\end{remark}
\begin{remark}\label{re1.2}
According to \eqref{B}, the upper bound of the temperature $\theta$ is not the key point to make sure that the solution $(\rho,\mathbf{u},P)$ is a global one, and it may go to infinity in the vacuum region within the life span of our strong solution.
\end{remark}
\begin{remark}\label{re1.3}
Compared with \cite{HX2016}, where the authors investigated a blow-up criterion for the IBVP of 3D full compressible Navier-Stokes equations with zero heat conduction, there is no need to impose additional restrictions on the viscosity coefficients $\mu$ and $\lambda$ except the physical restrictions \eqref{1.02}.
\end{remark}

\begin{remark}
As stated above, the same criterion has been derived by the author \cite{Z2017} for the Cauchy problem. However, the key idea in the present paper is different from that of in \cite{Z2017}. Roughly speaking, weighted energy estimates and Hardy-type inequalities play a crucial role in \cite{Z2017}, while the key ingredient of the analysis here is critical Sobolev inequalities of logarithmic type (see \eqref{2.2} and \eqref{zzz}).
\end{remark}

We now make some comments on the analysis of this paper. We mainly make use of continuation argument to prove Theorem \ref{thm1.1}. That is, suppose that \eqref{B} were false, i.e.,
\begin{equation*}
\lim_{T\rightarrow T^*}\left(\|\rho\|_{L^{\infty}(0,T;L^\infty)}
+\|P\|_{L^{\infty}(0,T;L^\infty)}\right)\leq M_0<\infty.
\end{equation*}
We want to show that
\begin{equation*}
\sup_{0\leq t\leq T^*}\left(\|(\rho,P)\|_{W^{1,q}}
+\|\nabla\mathbf{u}\|_{H^1}\right) \leq C<+\infty.
\end{equation*}
We first prove (see Lemma \ref{lem35}) that a control of upper bound of the density implies a control on the $L_t^\infty L_x^2$-norm of $\nabla\mathbf{u}$. To this end, the key ingredient of the analysis is a logarithmic Sobolev inequality (see Lemma \ref{lem25}). The inequality implies the uniform estimate of $\|\mathbf{u}\|_{L^2(0,T;L^\infty)}$ due to the a priori estimate of
$\|\mathbf{u}\|_{L^2(0,T;H^1)}$ from the energy estimate \eqref{3.2}. Then we obtain the key a priori estimates on $L^\infty(0,T;L^{q})$-norm of $(\nabla\rho,\nabla P)$ by solving a logarithmic Gronwall inequality based on a Br{\'e}zis-Waigner type inequality (see Lemma \ref{lem26}) and the a priori estimates we have just derived.

The rest of this paper is organized as follows. In Section \ref{sec2}, we collect some elementary facts and inequalities that will be used later. Section \ref{sec3} is devoted to the proof of Theorem \ref{thm1.1}.

\section{Preliminaries}\label{sec2}

In this section, we will recall some known facts and elementary inequalities that will be used frequently later.

We begin with the following Gagliardo-Nirenberg inequality (see \cite[Theorem 10.1, p. 27]{F2008}).
\begin{lemma}[Gagliardo-Nirenberg]\label{lem22}
Let $\Omega\subset\mathbb{R}^2$ be a bounded smooth domain.
Assume that $1\leq q,r\leq\infty$, and $j,m$ are arbitrary integers satisfying $0\leq j<m$. If $v\in W^{m,r}(\Omega)\cap L^q(\Omega)$, then we have
\begin{equation*}
\|D^jv\|_{L^p}\leq C\|v\|_{L^q}^{1-a}\|v\|_{W^{m,r}}^{a},
\end{equation*}
where
\begin{equation*}
-j+\frac{2}{p}=(1-a)\frac{2}{q}+a\left(-m+\frac{2}{r}\right),
\end{equation*}
and
\begin{equation*}
\begin{split}
a\in
\begin{cases}
[\frac{j}{m},1),\ \ \text{if}\ m-j-\frac{2}{r}\ \text{is an nonnegative integer},\\
[\frac{j}{m},1],\ \ \text{otherwise}.
\end{cases}
\end{split}
\end{equation*}
The constant $C$ depends only on $m,j,q,r,a$, and $\Omega$.
\end{lemma}

Next, we give some regularity results of the following Lam{\'e} system with Dirichlet boundary condition, the proof can be found in \cite[Proposition 2.1]{SWZ20111}.
\begin{align}\label{2.1}
\begin{cases}
\mu\Delta\mathbf{U}+(\mu+\lambda)\nabla\divv\mathbf{U}=\mathbf{F},\ \ x\in\Omega,\\
\mathbf{U}=\mathbf{0},\ \ x\in\partial\Omega.
\end{cases}
\end{align}
\begin{lemma}\label{lem23}
Let $q\geq2$ and $\mathbf{U}$ be a weak  solution of \eqref{2.1}. There exists a constant $C$ depending only on $q,\mu,\lambda$,
and $\Omega$ such that the following estimates hold:
\begin{itemize}
\item [$\bullet$] If $\mathbf{F}\in L^{q}(\Omega)$, then
\begin{equation*}
\|\mathbf{U}\|_{W^{2,q}}\leq C\|\mathbf{F}\|_{L^q};
\end{equation*}
\item [$\bullet$] If $\mathbf{F}\in W^{-1,q}(\Omega)$ (i.e., $\mathbf{F}= \divv\mathbf{f}$ with $\mathbf{f}=(f_{ij})_{3\times3}, f_{ij}\in L^q(\Omega)$), then
\begin{equation*}
\|\mathbf{U}\|_{W^{1,q}}\leq C\|\mathbf{f}\|_{L^q};
\end{equation*}
\item [$\bullet$] If $\mathbf{F}\in W^{-1,q}(\Omega)$ (i.e., $\mathbf{F}= \divv\mathbf{f}$ with $\mathbf{f}=(f_{ij})_{3\times3}, f_{ij}\in L^\infty(\Omega)$), then
\begin{equation*}
[\nabla\mathbf{U}]_{BMO}\leq C\|\mathbf{f}\|_{L^\infty}.
\end{equation*}
\end{itemize}
\end{lemma}
Here $BMO(\Omega)$ stands for the John-Nirenberg's space of bounded mean oscillation whose norm is defined by
\begin{align*}
\|f\|_{BMO}=\|f\|_{L^2}+[f]_{BMO}
\end{align*}
with
\begin{align*}
[f]_{BMO}=\sup_{x\in\Omega,r\in(0,d)}\frac{1}{|\Omega_r(x)|}
\int_{\Omega_r(x)}|f(y)-f_{\Omega_r}(x)|dy,
\end{align*}
and
\begin{align*}
f_{\Omega_r}(x)=\frac{1}{|\Omega_r(x)|}
\int_{\Omega_r(x)}f(y)dy,
\end{align*}
where $\Omega_r(x)=B_r(x)\cap\Omega$, $B_r(x)$ is the ball with center $x$
and radius $r$, and $d$ is the diameter of $\Omega$. $|\Omega_r(x)|$ denotes the Lebesgue measure of $\Omega_r(x)$.

Next, we decompose the velocity field into two parts in order to overcome the difficulty caused by the boundary, namely $\mathbf{u}=\mathbf{v}+\mathbf{w}$, where $\mathbf{v}$ is the solution to the Lam{\'e} system
\begin{align}\label{2.3}
\begin{cases}
\mu\Delta\mathbf{v}+(\mu+\lambda)\nabla\divv\mathbf{v}=\nabla P,\ \ x\in\Omega,\\
\mathbf{v}=\mathbf{0},\ \ x\in\partial\Omega,
\end{cases}
\end{align}
and $\mathbf{w}$ satisfies the following boundary value problem
\begin{align}\label{2.4}
\begin{cases}
\mu\Delta\mathbf{w}+(\mu+\lambda)\nabla\divv\mathbf{w}
=\rho\dot{\mathbf{u}},\ \ x\in\Omega,\\
\mathbf{w}=\mathbf{0},\ \ x\in\partial\Omega.
\end{cases}
\end{align}
By virtue of Lemma \ref{lem23}, one has the following key estimates for $\mathbf{v}$ and $\mathbf{w}$.
\begin{lemma}\label{lem21}
Let $\mathbf{v}$ and $\mathbf{w}$ be a solution of \eqref{2.3} and \eqref{2.4} respectively. Then for any $p\geq2$, there is a constant $C>0$ depending only on $p,\mu,\lambda$, and $\Omega$ such that
\begin{equation}\label{2.5}
\|\mathbf{v}\|_{W^{1,p}}\leq C\|P\|_{L^p},
\end{equation}
and
\begin{equation}\label{2.6}
\|\mathbf{w}\|_{W^{2,p}}
\leq C\|\rho\dot{\mathbf{u}}\|_{L^p}.
\end{equation}
\end{lemma}

Next, we state a critical Sobolev inequality of logarithmic type, which is originally
due to Br{\'e}zis-Wainger \cite{BW1980}.  The reader can refer to  \cite[Section 2]{HW20130} for the proof.
\begin{lemma}\label{lem25}
Assume $\Omega$ is a bounded smooth domain in $\mathbb{R}^2$ and $f\in L^2(s,t;W^{1,q}(\Omega))$ with some $q>2$ and $0\leq s<t\leq\infty$, then there is a constant $C>0$ depending only on $q$ and $\Omega$ such that
\begin{equation}\label{2.2}
\|f\|_{L^2(s,t;L^\infty)}^2\leq C\left(1+\|f\|_{L^2(s,t;H^1)}^2
\log(e+\|f\|_{L^2(s,t;W^{1,q})})\right).
\end{equation}
\end{lemma}

Finally, the following variant of the Br{\'e}zis-Waigner inequality
plays a crucial role in obtaining the estimate of $\|(\nabla\rho,\nabla P)\|_{L^q}$. For its proof, please refer to \cite[Lemma 2.3]{SWZ20111}.
\begin{lemma}\label{lem26}
Assume $\Omega$ is a bounded smooth domain in $\mathbb{R}^2$ and $f\in W^{1,q}(\Omega)$ with some $q>2$, then there is a constant $C>0$ depending only on $q$ and $\Omega$ such that
\begin{equation}\label{zzz}
\|f\|_{L^\infty}\leq C\left(1+\|f\|_{BMO}
\log(e+\|f\|_{W^{1,q}})\right).
\end{equation}
\end{lemma}

\section{Proof of Theorem \ref{thm1.1}}\label{sec3}

Let $(\rho,\mathbf{u},P)$ be a strong solution described in Theorem \ref{thm1.1}. Suppose that \eqref{B} were false, that is, there exists a constant $M_0>0$ such that
\begin{equation}\label{3.1}
\lim_{T\rightarrow T^*}\left(\|\rho\|_{L^{\infty}(0,T;L^\infty)}
+\|P\|_{L^{\infty}(0,T;L^\infty)}\right)\leq M_0<\infty.
\end{equation}

%

First, we have the following energy estimate.
\begin{lemma}\label{lem32}
Under the condition \eqref{3.1}, it holds that for any $T\in[0,T^*)$,
\begin{equation}\label{3.2}
\sup_{0\leq t\leq T}\|\sqrt{\rho}\mathbf{u}\|_{L^2}^2+\int_{0}^T\|\nabla\mathbf{u}\|_{L^2}^2 dt \leq C,
\end{equation}
where and in what follows, $C,C_1$ stand for generic positive constants depending only on $\Omega,M_0,\lambda,\mu,T^{*}$, and the initial data.
\end{lemma}
{\it Proof.}
1. It follows from \eqref{1.10}$_3$ that
\begin{equation}\label{3.3}
P_t+\mathbf{u}\cdot\nabla P+2P\divv\mathbf{u}
=F\triangleq2\mu|\mathfrak{D}(\mathbf{u})|^2
+\lambda(\divv\mathbf{u})^2\geq0.
\end{equation}
Define particle path before blow-up time
\begin{align*}
\begin{cases}
\frac{d}{dt}\mathbf{X}(x,t)
=\mathbf{u}(\mathbf{X}(x,t),t),\\
\mathbf{X}(x,0)=x.
\end{cases}
\end{align*}
Thus, along particle path, we obtain from \eqref{3.3} that
\begin{align*}
\frac{d}{dt}P(\mathbf{X}(x,t),t)
=-2P\divv\mathbf{u}+F,
\end{align*}
which implies
\begin{equation}\label{3.4}
P(\mathbf{X}(x,t),t)=\exp\left(-2\int_{0}^t\divv\mathbf{u}ds\right)
\left[P_0+\int_{0}^t\exp\left(2\int_{0}^s\divv\mathbf{u}d\tau\right)Fds\right]\geq0.
\end{equation}

2. Multiplying \eqref{1.10}$_2$ by $\mathbf{u}$ and integrating the resulting equation over $\Omega$, we deduce from Cauchy-Schwarz inequality that
\begin{align*}
\frac{1}{2}\frac{d}{dt}\int\rho|\mathbf{u}|^2dx
+\int\left[\mu|\nabla\mathbf{u}|^2+(\lambda+\mu)(\divv\mathbf{u})^2
\right]dx & =\int P\divv\mathbf{u}dx \notag \\
& \leq\frac{\mu}{2}\int|\nabla\mathbf{u}|^2dx
+\frac{1}{2\mu}\int P^2dx,
\end{align*}
which combined with  \eqref{3.1} implies that
\begin{align}\label{3.5}
\frac{d}{dt}\|\sqrt{\rho}\mathbf{u}\|_{L^2}^2
+\mu\|\nabla\mathbf{u}\|_{L^2}^2
\leq C.
\end{align}
So the desired \eqref{3.2} follows from \eqref{3.5} integrated with respect to $t$. This completes the proof of Lemma \ref{lem32}.
\hfill $\Box$

The following lemma gives the estimate on the spatial gradients of the velocity, which is crucial for deriving the higher order estimates of the solution.
\begin{lemma}\label{lem35}
Under the condition \eqref{3.1}, it holds that for any $T\in[0,T^*)$,
\begin{equation}\label{5.1}
\sup_{0\leq t\leq T}\|\nabla\mathbf{u}\|_{L^2}^{2}
+\int_{0}^{T}\|\sqrt{\rho}\dot{\mathbf{u}}\|_{L^2}^{2}
 dt \leq C.
\end{equation}
\end{lemma}
{\it Proof.}
1. Multiplying \eqref{1.10}$_2$ by $\mathbf{u}_t$ and integrating the resulting equation over $\Omega$ give rise to
\begin{align}\label{5.2}
\frac12\frac{d}{dt}\int\left(\mu|\nabla\mathbf{u}|^2
+(\mu+\lambda)(\divv\mathbf{u})^2\right)dx
+\int\rho|\dot{\mathbf{u}}|^2dx= \int\rho\dot{\mathbf{u}}\cdot\mathbf{u}\cdot\nabla\mathbf{u}dx
+\int P\divv\mathbf{u}_tdx.
\end{align}
It follows from Cauchy-Schwarz inequality and \eqref{3.1} that
\begin{align}\label{5.4}
\left|\int\rho\dot{\mathbf{u}}\cdot\mathbf{u}\cdot\nabla\mathbf{u}dx\right|
\leq \frac14\|\sqrt{\rho}\dot{\mathbf{u}}\|_{L^2}^2
+C\|\mathbf{u}\|_{L^\infty}^2\|\nabla\mathbf{u}\|_{L^2}^2.
\end{align}
To bound the second term of the right hand side of \eqref{5.2}, we decompose $\mathbf{u}$ into $\mathbf{u}=\mathbf{v}+\mathbf{w}$, where $\mathbf{v}$ and $\mathbf{w}$ satisfy \eqref{2.3} and \eqref{2.4}, respectively.
Then we obtain from \eqref{2.3} and integration by parts that
\begin{align}\label{5.3}
\int P\divv\mathbf{u}_tdx & = \frac{d}{dt}\left(\int P\divv\mathbf{u}dx\right)-\int P_t\divv\mathbf{u}dx
\notag \\
& = \frac{d}{dt}\left(\int P\divv\mathbf{u}dx\right)-\int P_t\divv\mathbf{v}dx
-\int P_t\divv\mathbf{w}dx
\notag \\
& = \frac{d}{dt}\left(\int P\divv\mathbf{u}dx\right)+\int \nabla P_t\cdot\mathbf{v}dx
-\int P_t\divv\mathbf{w}dx \notag \\
& = \frac{d}{dt}\left(\int P\divv\mathbf{u}dx\right)+\int (\mu\Delta\mathbf{v}+(\mu+\lambda)\nabla\mathbf{v})_t\cdot\mathbf{v}dx
-\int P_t\divv\mathbf{w}dx \notag \\
& = \frac12\frac{d}{dt}\int\left( 2P\divv\mathbf{u}-\mu|\nabla\mathbf{v}|^2-(\mu+\lambda)(\divv\mathbf{v})^2\right)dx
-\int P_t\divv\mathbf{w}dx.
\end{align}
Denote
\begin{equation*}
E\triangleq\theta+\frac12|\mathbf{u}|^2,
\end{equation*}
then it follows from \eqref{1.10} that $E$ satisfies
\begin{equation*}
(\rho E)_t+\divv(\rho\mathbf{u}E+P\mathbf{u})
=\frac12\mu\Delta|\mathbf{u}|^2+\mu\divv(\mathbf{u}\cdot\nabla\mathbf{u})
+\lambda\divv(\mathbf{u}\divv\mathbf{u}),
\end{equation*}
which combined with \eqref{2.6} and \eqref{3.1} gives that
\begin{align*}
-\int P_t\divv\mathbf{w}dx
& = -\int (\rho E)_t\divv\mathbf{w}dx
+\frac12\int (\rho|\mathbf{u}|^2)_t\divv\mathbf{w}dx
\notag \\
& =-\int \left(\rho E\mathbf{u}+P\mathbf{u}-\mu\nabla\mathbf{u}\cdot\mathbf{u}
-\mu\mathbf{u}\cdot\nabla\mathbf{u}-\lambda\mathbf{u}\divv\mathbf{u}\right)
\cdot\nabla\divv\mathbf{w}dx
\notag \\
& \quad -\frac12\int\divv(\rho\mathbf{u})|\mathbf{u}|^2\divv\mathbf{w}dx
+\int\rho\mathbf{u}\cdot\mathbf{u}_t\divv\mathbf{w}dx \notag \\
& = -\int \left(2P\mathbf{u}+\frac12\rho|\mathbf{u}|^2\mathbf{u}-\mu\nabla\mathbf{u}\cdot\mathbf{u}
-\mu\mathbf{u}\cdot\nabla\mathbf{u}-\lambda\mathbf{u}\divv\mathbf{u}\right)
\cdot\nabla\divv\mathbf{w}dx
\notag \\
& \quad +\frac12\int\rho|\mathbf{u}|^2\mathbf{u}\cdot\nabla\divv\mathbf{w}dx
+\int\rho\dot{\mathbf{u}}\cdot\mathbf{u}\divv\mathbf{w}dx\notag \\
& = -\int \left(2P\mathbf{u}-\mu\nabla\mathbf{u}\cdot\mathbf{u}
-\mu\mathbf{u}\cdot\nabla\mathbf{u}-\lambda\mathbf{u}\divv\mathbf{u}\right)
\cdot\nabla\divv\mathbf{w}dx
+\int\rho\dot{\mathbf{u}}\cdot\mathbf{u}\divv\mathbf{w}dx
\notag \\
& \leq C\|\mathbf{u}\|_{L^\infty}\left(\|P\|_{L^2}+\|\nabla\mathbf{u}\|_{L^2}\right)
\|\nabla^2\mathbf{w}\|_{L^2}+C\|\mathbf{u}\|_{L^\infty}\|\rho\dot{\mathbf{u}}\|_{L^2}
\|\nabla\mathbf{w}\|_{L^2}\notag \\
& \leq C\|\mathbf{u}\|_{L^\infty}\left(\|P\|_{L^2}+\|\nabla\mathbf{u}\|_{L^2}\right)
\|\rho\dot{\mathbf{u}}\|_{L^2}+C\|\mathbf{u}\|_{L^\infty}\|\rho\dot{\mathbf{u}}\|_{L^2}
\|\nabla\mathbf{w}\|_{L^2} \notag \\
& \leq \frac14\|\sqrt{\rho}\dot{\mathbf{u}}\|_{L^2}^2
+C\|\mathbf{u}\|_{L^\infty}^2\left(1+\|\nabla\mathbf{u}\|_{L^2}^2
+\|\nabla\mathbf{v}\|_{L^2}^2\right),
\end{align*}
which together with \eqref{5.3} yields
\begin{align}\label{5.30}
\int P\divv\mathbf{u}_tdx & \leq \frac12\frac{d}{dt}\int\left(2P\divv\mathbf{u}-\mu|\nabla\mathbf{v}|^2-(\mu+\lambda)
(\divv\mathbf{v})^2\right)dx \notag \\
& \quad +\frac14\|\sqrt{\rho}\dot{\mathbf{u}}\|_{L^2}^2
+C\|\mathbf{u}\|_{L^\infty}^2\left(1+\|\nabla\mathbf{u}\|_{L^2}^2
+\|\nabla\mathbf{v}\|_{L^2}^2\right).
\end{align}
Putting \eqref{5.4} and \eqref{5.30} into \eqref{5.2}, we get
\begin{align}\label{5.8}
B'(t)+\|\sqrt{\rho}\dot{\mathbf{u}}\|_{L^2}^2
\leq C\|\mathbf{u}\|_{L^\infty}^2\left(1+\|\nabla\mathbf{u}\|_{L^2}^2
+\|\nabla\mathbf{v}\|_{L^2}^2\right),
\end{align}
where
\begin{align}\label{5.9}
B(t) & \triangleq \frac{\mu}{2}\|\nabla\mathbf{u}\|_{L^2}^2
+\frac{\lambda+\mu}{2}\|\divv\mathbf{u}\|_{L^2}^2
+\frac{\mu}{2}\|\nabla\mathbf{v}\|_{L^2}^2
+\frac{\lambda+\mu}{2}\|\divv\mathbf{v}\|_{L^2}^2
-\int P\divv\mathbf{u}dx.
\end{align}

2. Let
\begin{align}\label{5.12}
\Phi(t)\triangleq e+\sup_{0\leq\tau\leq t}\left(\|\nabla\mathbf{u}\|_{L^2}^2+\|\nabla\mathbf{v}\|_{L^2}^2
\right)+\int_0^t\|\sqrt{\rho}\dot{\mathbf{u}}\|_{L^2}^{2}d\tau.
\end{align}
Then we obtain from Gronwall's inequality that for every $0\leq s\leq T<T^*$,
\begin{align}\label{5.13}
\Phi(T)\leq C\Phi(s)\exp\left\{C\int_s^T\|\mathbf{u}\|_{L^\infty}^2d\tau\right\}.
\end{align}
From Lemma \ref{lem25}, we get
\begin{align}\label{5.14}
\|\mathbf{u}\|_{L^2(s,T;L^\infty)}^2
& \leq C\left(1+\|\mathbf{u}\|_{L^2(s,T;H^1)}^2
\log\left(e+\|\mathbf{u}\|_{L^2(s,T;W^{1,3})}\right)\right) \notag \\
& \leq C_1\left(1+\|\nabla\mathbf{u}\|_{L^2(s,T;L^2)}^2
\log(C\Phi(T))\right),
\end{align}
where one has used the Poincar{\'e} inequality, \eqref{3.2}, and the following fact
\begin{align*}
\|\mathbf{u}\|_{W^{1,3}}^2\leq \|\mathbf{w}\|_{W^{1,3}}^2+\|\mathbf{v}\|_{W^{1,3}}^2 \leq C\|\mathbf{w}\|_{W^{2,2}}^2+C\|P\|_{L^{3}}^2\leq C\|\sqrt{\rho}\dot{\mathbf{u}}\|_{L^2}^{2}+C.
\end{align*}
The combination \eqref{5.13} and \eqref{5.14} gives rise to
\begin{align}\label{5.15}
\Phi(T)\leq C\Phi(s)(C\Phi(T))^{C_1\|\nabla\mathbf{u}\|_{L^2(s,T;L^2)}^2}.
\end{align}
Recalling \eqref{3.2}, one can choose $s$ close enough to $T^*$ such that
\begin{align*}
\lim_{T\rightarrow T^{*-}}C_1\|\nabla\mathbf{u}\|_{L^2(s,T;L^2)}^2\leq \frac12.
\end{align*}
Hence, for $s<T<T^*$, we have
\begin{align*}
\Phi(T)\leq C\Phi^2(s)<\infty.
\end{align*}
This completes the proof of Lemma \ref{lem35}.
\hfill $\Box$

Next, motivated by \cite{H1995}, we have the following estimates on the material derivatives of the velocity which are important for the higher order estimates of strong solutions.
\begin{lemma}\label{lem36}
Under the condition \eqref{3.1}, it holds that for any $T\in[0,T^*)$,
\begin{equation}\label{6.1}
\sup_{0\leq t\leq T}\|\sqrt{\rho}\dot{\mathbf{u}}\|_{L^2}^2
+\int_0^T\|\nabla\dot{\mathbf{u}}\|_{L^2}^2dt \leq C.
\end{equation}
\end{lemma}
{\it Proof.}
By the definition of $\dot{\mathbf{u}}$, we can rewrite $\eqref{1.10}_2$ as follows:
\begin{equation}\label{6.2}
\rho\dot{\mathbf{u}}+\nabla P
=\mu\Delta\mathbf{u}+(\lambda+\mu)\nabla\divv\mathbf{u}.
\end{equation}
Differentiating \eqref{6.2} with respect to $t$ and using \eqref{1.10}$_1$, we have
\begin{align}\label{6.3}
\rho\dot{\mathbf{u}}_t+\rho\mathbf{u}\cdot\nabla\dot{\mathbf{u}}
+\nabla P_t&=\mu\Delta\dot{\mathbf{u}}+(\lambda+\mu)\divv\dot{\mathbf{u}}
-\mu\Delta(\mathbf{u}\cdot\nabla\mathbf{u}) \nonumber \\
& \quad -(\lambda+\mu)\divv(\mathbf{u}\cdot\nabla\mathbf{u})
+\divv(\rho\dot{\mathbf{u}}\otimes\mathbf{u}).
\end{align}
Multiplying \eqref{6.3} by $\dot{\mathbf{u}}$ and integrating by parts over $\Omega$, we get
\begin{align}\label{6.4}
&\frac12\frac{d}{dt}\int \rho|\dot{\mathbf{u}}|^2\mbox{d}x+\mu\int|\nabla\dot{\mathbf{u}}|^2dx
+(\lambda+\mu)\int|\divv\dot{\mathbf{u}}|^2dx  \nonumber\\
&=\int \left(P_t\divv\dot{\mathbf{u}}+(\nabla P\otimes\mathbf{u}):\nabla\dot{\mathbf{u}}\right)dx
 +\mu\int [\divv(\Delta\mathbf{u}\otimes\mathbf{u})-\Delta(\mathbf{u}\cdot\nabla\mathbf{u})]\cdot \dot{\mathbf{u}}dx  \nonumber\\
& \quad+(\lambda+\mu)\int [(\nabla\divv\mathbf{u})\otimes\mathbf{u}
-\nabla\divv(\mathbf{u}\cdot\nabla\mathbf{u})]\cdot \dot{\mathbf{u}}dx\triangleq\sum_{i=1}^{3}J_i,
\end{align}
where $J_i$ can be bounded as follows.

It follows from $\eqref{1.10}_3$ that
\begin{align}\label{6.5}
J_1 & =\int \left(-\divv(P\mathbf{u})\divv\dot{\mathbf{u}}
-P\divv\mathbf{u}\divv\dot{\mathbf{u}}
+\mathcal{T}(\mathbf{u}):\nabla\mathbf{u}\divv\dot{\mathbf{u}}
+(\nabla P\otimes\mathbf{u}):\nabla\dot{\mathbf{u}}\right)dx   \nonumber\\
&=\int \left(P\mathbf{u}\nabla\divv\dot{\mathbf{u}}-P\divv\mathbf{u}\divv\dot{\mathbf{u}}
+\mathcal{T}(\mathbf{u}):\nabla\mathbf{u}\divv\dot{\mathbf{u}}
-P\nabla\mathbf{u}^\top:\nabla\dot{\mathbf{u}}
-P\mathbf{u}\nabla\divv\dot{\mathbf{u}}\right)dx \nonumber\\
&=\int \left(-P\divv\mathbf{u}\divv\dot{\mathbf{u}}
+\mathcal{T}(\mathbf{u}):\nabla\mathbf{u}\divv\dot{\mathbf{u}}
-P\nabla\mathbf{u}^\top:\nabla\dot{\mathbf{u}}\right)dx\nonumber\\
&\leq C\int \left(|\nabla\mathbf{u}||\nabla\dot{\mathbf{u}}|
+|\nabla\mathbf{u}|^2|\nabla\dot{\mathbf{u}}|\right)dx  \nonumber\\
&\leq C\left(\|\nabla\mathbf{u}\|_{L^2}+\|\nabla\mathbf{u}\|_{L^4}^2\right)
\|\nabla\dot{\mathbf{u}}\|_{L^2},
\end{align}
where $\mathcal{T}(\mathbf{u})=2\mu\mathfrak{D}(\mathbf{u})+\lambda\divv\mathbf{u}\mathbb{I}_2$.

For $J_2$ and $J_3$, notice that for all $1\leq i,j,k\leq 2,$ one has
\begin{align}
\partial_j(\partial_{kk}u_iu_j)-\partial_{kk}(u_j\partial_j u_i)&=\partial_k(\partial_ju_j\partial_ku_i)-\partial_k(\partial_ku_j\partial_ju_i)
-\partial_j(\partial_ku_j\partial_ku_i),\nonumber\\
\partial_j(\partial_{ik}u_ku_j)-\partial_{ij}(u_k\partial_ku_j)
&=\partial_i(\partial_ju_j\partial_ku_k)-\partial_i(\partial_ju_k\partial_ku_j)
-\partial_k(\partial_iu_k\partial_ju_j).\nonumber
\end{align}
So integrating by parts gives
\begin{align}
J_2&=\mu\int [\partial_k(\partial_ju_j\partial_ku_i)-\partial_k(\partial_ku_j\partial_ju_i)
-\partial_j(\partial_ku_j\partial_ku_i)]\dot{u_i}\mbox{d}x\nonumber\\
&\leq C\|\nabla\mathbf{u}\|_{L^4}^2\|\nabla\dot{\mathbf{u}}\|_{L^2},\label{6.7}\\
J_3&=(\lambda+\mu)\int [\partial_i(\partial_ju_j\partial_ku_k)-\partial_i(\partial_ju_k\partial_ku_j)
-\partial_k(\partial_iu_k\partial_ju_j)]\dot{u_i}\mbox{d}x\nonumber\\
&\leq C\|\nabla\mathbf{u}\|_{L^4}^2\|\nabla\dot{\mathbf{u}}\|_{L^2}.\label{6.8}
\end{align}
Inserting \eqref{6.5}--\eqref{6.8} into \eqref{6.4} and applying \eqref{5.1} lead to
\begin{align}\label{6.9}
& \frac12\frac{d}{dt}
\|\sqrt{\rho}\dot{\mathbf{u}}\|_{L^2}^2+\mu\|\nabla\dot{\mathbf{u}}\|_{L^2}^2
+(\lambda+\mu)\|\divv\dot{\mathbf{u}}\|_{L^2}^2  \nonumber\\
& \leq C(\|\nabla\mathbf{u}\|_{L^2}+\|\nabla\mathbf{u}\|_{L^4}^2)
\|\nabla\dot{\mathbf{u}}\|_{L^2}\nonumber\\
&\leq \frac{\mu}{2}\|\nabla\dot{\mathbf{u}}\|_{L^2}^2 +C(\mu)\left(\|\nabla\mathbf{u}\|_{L^4}^4+1\right),
\end{align}
which implies
\begin{align}\label{6.16}
\frac{d}{dt}\|\sqrt{\rho}\dot{\mathbf{u}}\|_{L^2}^2
+\mu\|\nabla\dot{\mathbf{u}}\|_{L^2}^2
\leq C\|\nabla\mathbf{u}\|_{L^4}^4+C.
\end{align}
By virtue of Gagliardo-Nirenberg inequality, \eqref{2.5}, \eqref{2.6}, and \eqref{5.1}, one has
\begin{align}\label{6.17}
\|\nabla\mathbf{u}\|_{L^4}^4
& \leq C\|\nabla\mathbf{v}\|_{L^4}^4+C\|\nabla\mathbf{w}\|_{L^4}^4 \notag \\
& \leq C\|P\|_{L^4}^4+C\|\nabla\mathbf{w}\|_{L^2}^2\|\nabla\mathbf{w}\|_{H^1}^2 \notag \\ & \leq C
+C\|\rho\dot{\mathbf{u}}\|_{L^2}^2 \notag \\ &
\leq C+C\|\sqrt{\rho}\dot{\mathbf{u}}\|_{L^2}^2.
\end{align}
Consequently, we obtain the desired \eqref{6.1} from \eqref{6.16}, \eqref{6.17}, and Gronwall's inequality. This completes the proof of Lemma \ref{lem36}.
\hfill $\Box$

The following lemma will treat the higher order derivatives of the solutions which are needed to guarantee the extension of local strong solution to be a global one.
\begin{lemma}\label{lem37}
Under the condition \eqref{3.1}, and let $q>2$ be as in Theorem \ref{thm1.1}, then it holds that for any $T\in[0,T^*)$,
\begin{equation}\label{7.1}
\sup_{0\leq t\leq T}\left(\|(\rho,P)\|_{W^{1,q}}+\|\nabla\mathbf{u}\|_{H^1}\right)\leq C.
\end{equation}
\end{lemma}
{\it Proof.}
1. For $q>2$, it follows from the mass equation \eqref{1.10}$_1$ that $\nabla\rho$ satisfies
\begin{align}\label{7.2}
\frac{d}{dt}\|\nabla\rho\|_{L^q}
& \leq C(q)(1+\|\nabla\mathbf{u} \|_{L^\infty})\|\nabla\rho\|_{L^q}
+C(q)\|\nabla^2\mathbf{u}\|_{L^q} \nonumber \\
& \leq C(1+\|\nabla\mathbf{w}\|_{L^\infty}+\|\nabla\mathbf{v}\|_{L^\infty})
\|\nabla\rho\|_{L^q}
+C\left(\|\nabla^2\mathbf{w}\|_{L^q}+\|\nabla^2\mathbf{v}\|_{L^q}\right) \nonumber \\
& \leq C(1+\|\nabla\mathbf{w}\|_{L^\infty}+\|\nabla\mathbf{v}\|_{L^\infty})
\|\nabla\rho\|_{L^q}
+C\|\nabla^2\mathbf{w}\|_{L^q}+C\|\nabla P\|_{L^q}
\end{align}
due to the following fact
\begin{equation*}
\|\nabla^2\mathbf{v}\|_{L^q}
\leq C\|\nabla P\|_{L^q},
\end{equation*}
which follows from the standard $L^q$-estimate for the following elliptic system
\begin{align*}
\begin{cases}
\mu\Delta\mathbf{v}+(\lambda+\mu)\nabla\divv\mathbf{v}
=\nabla P,\ \ x\in\Omega,\\
\mathbf{v}=\mathbf{0},\  x\in\partial\Omega.
\end{cases}
\end{align*}
From Lemma \ref{lem26} and \eqref{2.5}, one gets
\begin{align}\label{7.3}
\|\nabla\mathbf{v}\|_{L^\infty}
& \leq C\left(1+\|\nabla\mathbf{v}\|_{BMO}\log(e+\|\nabla\mathbf{v}\|_{W^{1,q}})\right) \nonumber \\
& \leq C\left(1+(\|P\|_{L^2}+\|P\|_{L^\infty})
\log(e+\|\nabla\mathbf{v}\|_{W^{1,q}})\right)\nonumber \\
& \leq C\left(1+\log(e+\|\nabla P\|_{L^{q}})\right).
\end{align}
By virtue of Sobolev's embedding theorem, \eqref{2.6}, and \eqref{3.1}, one deduces that
\begin{align}\label{7.5}
\|\nabla\mathbf{w}\|_{L^\infty}
\leq \|\mathbf{w}\|_{W^{2,q}}\leq C\|\rho\dot{\mathbf{u}}\|_{L^q}
\leq C\|\dot{\mathbf{u}}\|_{L^q}
\leq C\|\nabla\dot{\mathbf{u}}\|_{L^2}.
\end{align}
Moreover, we have
\begin{align}\label{7.6}
\|\nabla^2\mathbf{w}\|_{L^q}
\leq \|\mathbf{w}\|_{W^{2,q}}\leq C\|\nabla\dot{\mathbf{u}}\|_{L^2}.
\end{align}
Substituting \eqref{7.3}--\eqref{7.6} into \eqref{7.2}, we derive that
\begin{align}\label{7.4}
\frac{d}{dt}\|\nabla\rho\|_{L^{q}}\leq C\left(1+\|\nabla\dot{\mathbf{u}}\|_{L^2}+\log(e+\|\nabla P\|_{L^{q}})\right)
\|\nabla\rho\|_{L^{q}}+C\|\nabla\dot{\mathbf{u}}\|_{L^2}+C\|\nabla P\|_{L^q}.
\end{align}
Similarly, one deduces from \eqref{1.10}$_3$ that $\nabla P$ satisfies for any $q>2$,
\begin{align}\label{7.7}
\frac{d}{dt}\|\nabla P\|_{L^q} & \leq C(q)(1+\|\nabla\mathbf{u}\|_{L^{\infty}})(\|\nabla P\|_{L^q}+\|\nabla^2\mathbf{u}\|_{L^q}) \notag \\
& \leq C\left(1+\|\nabla\dot{\mathbf{u}}\|_{L^2}+\log(e+\|\nabla P\|_{L^{q}})\right)
\|\nabla P\|_{L^{q}}+C\|\nabla\dot{\mathbf{u}}\|_{L^2}+C\|\nabla P\|_{L^q}.
\end{align}

2. Let
\begin{align*}
f(t)&\triangleq e+\|\nabla\rho\|_{L^{ q}}+\|\nabla P\|_{L^{ q}},\\
g(t)&\triangleq (1+\|\nabla\dot{\mathbf{u}}\|_{L^2})
\log(e+\|\nabla\dot{\mathbf{u}}\|_{L^2}).
\end{align*}
Then we derive from \eqref{7.4} and \eqref{7.7} that
\begin{equation}\label{7.11}
f'(t)\leq  Cg(t)f(t)\log{f(t)}+Cg(t)f(t)+Cg(t),
\end{equation}
which yields
\begin{equation}\label{7.12}
(\log f(t))'\leq Cg(t)+Cg(t)\log f(t)
\end{equation}
due to $f(t)>1.$
Thus it follows from  \eqref{7.12}, \eqref{6.1}, and  Gronwall's inequality that
\begin{equation}\label{7.13}
\sup_{0\leq t\leq T}\|(\nabla\rho,\nabla P)\|_{L^{q}}\leq C.
\end{equation}
We infer from \eqref{2.5}, \eqref{2.6}, and \eqref{3.1} that
\begin{align*}
\|\nabla^2\mathbf{u}\|_{L^2}
\leq \|\nabla^2\mathbf{v}\|_{L^2}+\|\nabla^2\mathbf{w}\|_{L^2}
\leq C\|\nabla P\|_{L^2}+C\|\sqrt{\rho}\dot{\mathbf{u}}\|_{L^2},
\end{align*}
which combined with \eqref{7.13}, \eqref{6.1}, and H{\"o}lder's inequality implies that
\begin{equation}\label{7.16}
\sup_{0\leq t\leq T}\|\nabla^2\mathbf{u}\|_{L^2}\leq C.
\end{equation}
Thus the desired \eqref{7.1} follows from  \eqref{7.13}, \eqref{7.16}, \eqref{5.1}, and \eqref{3.1}. The proof of Lemma \ref{lem37} is finished.  \hfill $\Box$

With Lemmas \ref{lem32}--\ref{lem37} at hand, we are now in a position to prove Theorem \ref{thm1.1}.

\textbf{Proof of Theorem \ref{thm1.1}.}
We argue by contradiction. Suppose that \eqref{B} were false, that is, \eqref{3.1} holds. Note that the general constant $C$ in Lemmas \ref{lem32}--\ref{lem37} is independent of $t<T^{*}$, that is, all the a priori estimates obtained in Lemmas \ref{lem32}--\ref{lem37} are uniformly bounded for any $t<T^{*}$. Hence, the function
\begin{equation*}
(\rho,\mathbf{u},P)(x,T^{*})
\triangleq\lim_{t\rightarrow T^{*}}(\rho,\mathbf{u},P)(x,t)
\end{equation*}
satisfy the initial condition \eqref{2.02} at $t=T^{*}$.

Furthermore, standard arguments yield that $\rho\dot{\mathbf{u}}\in C([0,T];L^2)$, which
implies $$ \rho\dot{\mathbf{u}}(x,T^\ast)=\lim_{t\rightarrow
T^\ast}\rho\dot{\mathbf{u}}\in L^2. $$
Hence, $$-\mu\Delta{\mathbf{u}}-(\lambda+\mu)\nabla\mbox{div}\mathbf{u}+\nabla P
|_{t=T^\ast}=\sqrt{\rho}(x,T^\ast)g(x)
$$ with $$g(x)\triangleq
\begin{cases}
\rho^{-1/2}(x,T^\ast)(\rho\dot{\mathbf{u}})(x,T^\ast),&
\mbox{for}~~x\in\{x|\rho(x,T^\ast)>0\},\\
0,&\mbox{for}~~x\in\{x|\rho(x,T^\ast)=0\},
\end{cases}
$$
satisfying $g\in L^2$ due to \eqref{7.1}.
Therefore, one can take $(\rho,\mathbf{u},P)(x,T^\ast)$ as
the initial data and extend the local
strong solution beyond $T^\ast$. This contradicts the assumption on
$T^{\ast}$. Thus we finish the proof of Theorem \ref{thm1.1}.
\hfill $\Box$


\end{document}